\documentclass[11pt,reqno]{amsart}
\usepackage[dvips]{epsfig}
\RequirePackage{amsmath}
\usepackage{amscd,amsthm,amsfonts,amsopn,amssymb,verbatim}
\numberwithin{equation}{section}
\newtheorem{theorem}{Theorem} 

\newtheorem{lemma}[theorem]{Lemma}

\theoremstyle{remark}

\DeclareMathOperator{\supp}{supp\,}

\def\be{\begin{equation}}
\def\ee{\end{equation}}
\def\vp{\varphi}
\def\ve{\varepsilon}

\allowdisplaybreaks
\def\vp{\varphi}
\def\ve{\varepsilon}

\def\mod{\text{mod\,}}

\allowdisplaybreaks

\begin{document}
\begin{title}
{Decoupling inequalities and some mean-value theorems}
\end{title}
\author{Jean Bourgain}
\address{School of Mathematics,
Institute for Advanced Study, 
1 Einstein Drive, Princeton, NJ 08540}

\begin{abstract}
The purpose of this paper is to present some further applications of the general decoupling theory from [B-D1, 2] to
certain diophantine issues.  In particular, we consider mean value estimates relevant to the Bombieri-Iwaniec approach
to exponential sums and arising in the work of Robert and Sargos \cite{R-S}.
Our main input is a new mean value theorem. 
\end{abstract}
\maketitle

\centerline
{\scshape 0. Summary}

The aim of this Note is to illustrate how a version of the general decoupling inequality for hypersurfaces established in
\cite {B-D} permits to recover certain known mean-value theorems in number theory and establish some new ones.
Easy applications in this direction were already pointed out in \cite{B-D} and the material presented here is a further
development.
Our main emphasis will be on the method rather than the best exponents that can be obtained this way.

In the first section, we state a form  of the main decoupling theorem from \cite{B-D} to the situation of smooth hyper
surfaces in $\mathbb R^n$ with non-degenerate (but not necessarily definite) second fundamental form (a detailed argument
appears in \cite{B-D2}).
The motivation for this appears in Sections 2 and 3,  which aims at proving decoupling inequalities for real analytic curves
$\Gamma\subset\mathbb R^n$ not contained in a hyperplane.
The assumption of real analyticity is purely for convenience (it suffices for the subsequent applications) and a similar
result also holds in the smooth category.
An $(n-1)$-fold convolution of $\Gamma$ leads indeed to a hypersurface $S\subset\mathbb R^n$ of non-vanishing curvature.
The relevant statement is inequality \eqref{3.2} below with moment $q=2(n+1)$, where we consider the multi-linear
(i.e. $(n-1)$-linear) setting.
The next step is to reformulate this inequality as a  mean-value theorem for exponential sums stated as Theorem 1,
which is a quite general and optimal result.
Our first focus point are certain mean value inequalities arising
in the Bombieri-Iwaniec approach [B-I1, 2] to exponential sums and the subsequent
developments of this technique (see \cite{H} for the complete exposition).
More specifically, Theorem 1 is relevant to the so-called `first spacing problem' which is analytically captured by mean-value
expressions of the type
$$
N_8(\delta) =\int_0^1\int_0^1\int_0^1 \Big|\sum_{n\sim N} e(x_0 n+x_1n^2+x_2\frac 1\delta \Big(\frac
nN\Big)^{\frac 32}\Big)\Big|^8dx_0dx_1dx_2\eqno{(0.1)}
$$
$$
N_{10} (\delta, N\delta)=\int_0^1\int_0^1\int_0^1\int_0^1 \Big|\sum_{n\sim N}e(x_0n+x_1n^2+x_2\frac 1\delta\Big(\frac
nN\Big)^{\frac 32}+x_3 \frac 1{N\delta}\Big(\frac nN\Big)^{\frac 12}\Big|^{10}dx_0dx_1dx_2dx_3\eqno{(0.2)}
$$
and
$$
N_{12}(\delta, N\delta)=\int_0^1\int_0^1\int_0^1\int_0^1 \Big|\sum_{n\sim N} e(x_0n+x_1n^2+x_2\frac 1\delta\Big(\frac
nN\Big)^{3/2} +x_3 \frac 1{N\delta} \Big(\frac nN\big)^{\frac 12}\Big|^{12} dx_0dx_1dx_2 dx_3.\eqno{(0.3)}
$$
In the application, the most important range of $\delta$ is $\delta \sim \frac 1{N^2}$.
As a special case of a more general result, it was proven in \cite{B-I2} that
$$
N_8(\delta) \ll \delta N^{5+\ve}+N^{4+\ve}\eqno{(0.4)}
$$
and in \cite{H-K} that
$$
N_{10} (\delta, N\delta)\ll \delta N^{7+\ve}+ N^{5+\ve}.\eqno{(0.5)}
$$

Our Fourier analytical approach gives a quite different treatment and unified approach  to this problem.
In particular, Theorem 9 in Section 4 below shows that in fact
$$
N_{10}(\delta, N\delta)\ll N^{5+\ve} \ \text { for } \ \delta<N^{-\frac {33}{18}}.
\eqno{(0.6)}
$$
Since however the main contribution (at least in the treatment \cite{H}) in the exponential sum problem 
$$
\sum_{m\sim M} e\Big(TF\Big(\frac mM\Big)\Big)\eqno{(0.7)}
$$
has $\delta =\frac 1{N^2}$, the improvement (0.6) does not lead to new results on this matter.

Our next application are certain mean value results in the work of Robert and Sargos
\cite{R-S}.
It is proven in \cite{R-S} that
$$
I_6(N^{-3} ) =\int^1_0\int^1_0 \Big|\sum_{n\sim N} e(n^2 x+N^{-3} n^4 y)\Big|^6 dxdy\ll N^{3+\ve}
\eqno{(0.11)}
$$
$$
I_8(N^{-\frac 52} ) =\int^1_0\int^1_0 \Big(\sum_{n\sim N} e(n^2x+N^{-\frac 52} n^4 y)\Big|^8 dxdy\ll N^{\frac 92+\ve}
\eqno{(0.12)}
$$
$$
I_{10} (N^{-\frac {17}{8}})=\int^1_0\int^1_0 \Big|\sum_{n\sim N} e(n^2 x+N^{-\frac {17} 8} n^4 y)\big|^{10} dxdy\ll N^{\frac
{49}{8}+\ve}
\eqno{(0.13)}
$$
Inequality (0.11) is the optimal statement for the 6th moment (a different proof using the decoupling theorem for curves
appears in \cite{B-D}).
While (0.12), (0.13) are essentially sharp, they are not the optimal results for the 8th and 10th moment
respectively.
Since $I_p(\lambda)$ is a decreasing function of $\lambda$ for $p$ an even integer, (0.13) obviously implies that
$$
I_{10} =\int^1_0\int^1_0 \Big|\sum_{n\sim N} e(n^2x+ n^4y)\big|^{10} dxdy \ll N^{\frac {49}8+\ve}.
\eqno{(0.14)}
$$
In \cite{R-S} an application of (0.14) to Weyl's inequity is given, following a method initiated by Heath-Brown.
In view of the present state of the art, the relevant statement is the bound
$$
|f_8 (\alpha; N)|= \Big|\sum_{1\leq n\leq N} e(\alpha n^8)\Big| \ll N^{1-3.2^{-8}}(N^4 q^{-1} +1+qN^{-4})^{\frac 1{160}}
\eqno{(0.15)}
$$
assuming $|\alpha -\frac aq|\leq q^{-2}, q\geq 1, (a, q)=1$ (though the exponent $\sigma(8) = 3.2^{-8}=0,01171\cdots$ is superseded
by a recent result of Wooley, see Theorem 7.3 in \cite{W2}, which gives in particular $\sigma(8) =\frac 1{2.7.6}=\frac 1{84} =
0,01190\cdots$).

More recently, inequality (0.14) has been improved in \cite{P} to
$$
I_{10}\ll N^{6+\ve}
\eqno{(0.16)}
$$
using a different more arithmetical approach.
As a consequence the first factor in the r.h.s. of (0.15) is replaced by $N^{1-\frac {16}{5}.2^{-8}}$, i.e. $\sigma(8)=\frac 1{18}=
0,0125\cdots$.

In the final section of this paper, we establish the bounds
$$
I_8(N^{-\frac 73})\ll N^{\frac {13}3+\ve}\eqno{(0.17)}
$$
$$
I_{10} \leq I_{10} (N^{-\frac {5}3})\ll N^ {\frac{17} 3+\ve}
\eqno{(0.18)}
$$
implying a corresponding improvement $\sigma (8) =\frac {56}{15} 2^{-8} =0,0145\ldots$ in Weyl's inequality.

\section
{Decoupling inequality for smooth hypersurfaces with non-vanishing curvature}

Let us start by recalling the main result from \cite{B-D}, which is the so-called $\ell^2$-decoupling theorem for the Fourier transform of
distributions carried by hypersurfaces in $\mathbb R^n$ of positive curvature.
This is a quite general harmonic analysis result with diverse applications, in particular to PDE's and spectral theory (see \cite{B-D}
for some of these).

In order to formulate the result, we need some terminology.
Let $S\subset\mathbb R^n$ be a compact smooth hypersurface of positive curvature and denote $S_\delta$ ($\delta>0$ a small
parameter) a $\delta$-neighborhood of $S$.
Decompose $S_\delta$ as a union of tangent $\underbrace {\sqrt\delta\times\cdots\times \sqrt\delta}_{n-1}\times \delta$ boxes $\tau$
with bounded overlap.
Denoting $B_R\subset\mathbb R^n$ a ball of radius $R$, the following inequality holds for functions $f$ s.t. supp\,$\hat f\subset
S_\delta$
\be\label{1.1}
\Vert f\Vert_{L^p(B_{\frac 1\delta})} \ll \delta^{-\ve}\Big(\sum_\tau \Vert f_\tau\Vert^2_{L^p(B_{\frac 1\delta})}\Big)^{\frac 12}
\ \text { with } p=\frac {2(n+1)}{n-1}
\ee
and $f_\tau =(\hat f|_\tau)^\vee$ denoting the Fourier restriction of $f$ to the tile $\tau$.

By interpolation, \eqref{1.1} of course also holds for $2\leq p\leq \frac {2(n+1)}{n-1}$ while for $\frac {2(n+1)}{n-1}\leq p\leq
\infty$, the inequality becomes
\be\label{1.2}
\Vert f\Vert_{L^p(B_{\frac 1\delta})} \ll \delta^{-\frac {n-1}4+\frac {n+1}{2p}-\ve} \Big(\sum_\tau \Vert f_\tau\Vert^2_{L^p(B_{\frac
1\delta})}\Big)^{\frac 12}.
\ee
Next, let us relax the assumption on $S$, requiring $S$ to have non-degenerate (but not necessarily definite) second fundamental
form.
A statement such as \eqref{1.1} can not be valid any more.
For instance, if $S\subset\mathbb R^3$ is a ruled surface, we may take supp\,$\hat f$ in a $\sqrt\delta$-neighborhood of a straight
line segment with only the obvious decoupling available.
This problem of curvature break-down for lower dimensional sections of $S$ can be bypassed by a suitable reformulation of the 
decoupling  property.
Assuming $S$ as above and supp\,$\hat f\subset S_\delta$, one has for $\frac {2(n+1)}{n-1}\leq p\leq\infty$
\be\label{1.3}
\Vert f\Vert_{L^p(B_{\frac 1\delta})} \ll \delta^{-\frac {n-1}2+\frac np-\ve} \big(\sum_\tau \Vert f_\tau\Vert^p_{L^p(B_{\frac
1\delta})}\Big)^{\frac 1p}.
\ee
This statement is weaker than \eqref{1.2} but will perform equally well in what follows because in the applications below 
supp\,$\hat f$ will be uniformly spread out over $S$.

The proof of \eqref{1.3} requires a modification of the argument in \cite{B-D} (for positive curvature).
Details appear in \cite{B-D2}.
Our next goal is to derive from \eqref{1.3} a decoupling inequality for
curves $\Gamma\subset\mathbb R^n$ not lying in a hyperplane and which will imply our Theorem 1.

\section
{Construction of hypersurfaces from curves}

Let $\Gamma\subset \mathbb R^n$ be parametrized by $\Phi:[0, 1]\to\mathbb R^n: t\to \big(t, \vp_1(t), \ldots,
\vp_{n-1}(t)\big)$ where we assume for simplicity that $\vp_1, \ldots, \vp_{n-1}$ are real analytic and
(importantly) that $1, t, \vp_1, \ldots, \vp_{n-1}$ linearly independent.
In particular, $\Gamma$ does not lie in a hypersurface.
Our assumption means non-vanishing of the Wronskian determinant
\be\label{2.1}
W(\vp_1'', \ldots, \vp_{n-1}'')\not= 0.
\ee
We build a hypersurface $S\subset\mathbb R^n$ as $(n-1)$-fold sum set
\be\label{2.2}
S=\Gamma_1+\cdots+\Gamma_{n-1}
\ee
where $\Gamma_j= \Phi(I_j)$ and $I_1, \ldots, I_{n-1} \subset I\subset [0, 1]$ are fixed consecutive disjoint
subintervals.
Hence $S$ is parametrized by
\be\label{2.3}
\begin{cases}
x_0=t_1+\ldots t_{n-1}\\
x_1=\phi_1(t_1)+\cdots+\phi_1(t_{n-1})\\
\ \vdots\\
x_{n-1} =\phi_{n-1} (t_1)+\cdots+ \phi_{n-1} (t_{n-1})
\end{cases}
\ee
with $t_j\in I_j$.
Our aim is to show that the second fundamental form of $S$ is non-degenerate (but note that it may be
indefinite).
\bigskip

Perturb $t=(t_1, \ldots, t_{n-1}) \in I_1\times\cdots\times I_{n-1}$ to $(t_1+s_1, \ldots, t_{n-1} +s_{n-1})$,
$|s_j| =o(1)$, obtaining
{\small
\be\label{2.4}
\begin{cases}
x_0-t_1-\cdots-t_{n-1} \equiv x_0' = s_1+\cdots+s_{n-1}\\
{}\\
\begin{pmatrix}
x_1-\phi_1(t_1) -\cdots - \phi_1(t_{n-1})\\
\vdots\\
x_{n-1}-\phi_{n-1} (t_1) -\cdots -\phi_{n-1} (t_{n-1})\end{pmatrix}
\equiv 
\begin{pmatrix} x_1'\\ \vdots\\ x_{n-1}'\end{pmatrix} = D_1
\begin{pmatrix} s_1\\\vdots\\ s_{n-1}\end{pmatrix} +\frac 12 D_2 \begin{pmatrix} s_1^2\\ \vdots\\ s_{n-1}^2\end{pmatrix}
+ O(|s|^3)
\end{cases}
\ee}

with
\be\label{2.5}
D_1=
\left[\begin{matrix} \phi_1'(t_1)\cdots \phi_1'(t_{n-1})\\ \vdots\qquad\qquad \vdots\\ \phi_{n-1}'(t_1)\cdots \phi_{n-1}'
(t_{n-1})\end{matrix}\right]
\text { \ and \ } D_2=\left[
\begin{matrix} \phi_1''(t_1) \cdots \phi_1''(t_{n-1})\\ \vdots\qquad\qquad\qquad \\ \phi_{n-1} ''' (t_1) \cdots \phi_{n-1}'' (t_{n-1})\end{matrix}\right].
\ee
The non-vanishing of det$D_1$ can be derived from the non-vanishing of \hfill\break
$W(\phi_1', \ldots, \phi_{n-1}')$ which is a consequence of our
assumption \eqref{2.1}.

Hence, since $D_1$ is invertible and denoting $\xi=(1, \ldots, 1) \in\mathbb R^{n-1}$, the first equation in \eqref{2.4} gives
\begin{align}\label{2.6}
x_0'&= \Big\langle D_1^{-1} \begin{pmatrix} x_1'\\ \vdots\\ x_{n-1}'\end{pmatrix}, \xi \Big\rangle -\frac 12 \sum_{j=1}^{n-1} s_j^2\langle
D_1^{-1} D_2 e_j, \xi\rangle +O(|s|^3)\nonumber\\
&= \Big\langle \begin{pmatrix} x_1''\\ \vdots\\ x_{n-1}''\end{pmatrix}, \xi \Big\rangle -\frac 12\sum^{n-1}_{j=1} (x_j'')^2 \langle
D_1^{-1} D_2 e_j, \xi\rangle +O(|x''|^3)
\end{align}
where
$$
\begin{pmatrix} x_1''\\ \vdots\\ x_{n-1}''\end{pmatrix} =D_1^{-1} \begin{pmatrix} x_1'\\ \vdots\\
x_{n-1}'\end{pmatrix}.
$$
From \eqref{2.6}, it remains to ensure that
\be\label{2.7}
\langle D_1^{-1} D_2 e_j, \xi\rangle \not= 0 \text { for each } j=1, \ldots, n-1.
\ee
Take $j=1$.
Up to a multiplicative factor,

\be\label{2.8}
\langle D_2 e_1, (D_1^{-1})^* \xi\rangle \dot = \sum_{k=1}^{n-1} \phi_k''(t_1) \left|
\begin{matrix}
\phi_1'(t_1)\cdots\phi_1' (t_{n-1})\\ 
 \operatornamewithlimits 1\limits^{\vdots}_{\vdots} \qquad \cdots \qquad 1\\ 
\phi_{n-1}'(t_1)\cdots \phi_{n-1}'(t_{n-1})\end{matrix}
\right| \!\!\!\!\! \leftarrow k.
\ee
By the mean-value theorem, we obtain separated $t_1<t_2'< \cdots < t_{n-1}'$ such that
$$
\eqref{2.8} =\sum^{n-1}_{k=1} (-1)^k \phi_k'' (t_1) \left|\begin{matrix}
\phi_1'' (t_2')\ldots \phi_1'' (t_{n-1}')\\ \cdots \cdots\cdots\cdots\cdots\\ 
\phi_{n-1}'' (t_2') \ldots \phi_{n-1}'' (t_{n-1}')\end{matrix}\right|
\!\!\!\!\! \leftarrow k
$$
\vskip .3 true in
$$
= \left|\begin{matrix} &\phi_1''(t_1)  &\phi_1''(t_2') \ \ \cdots \!\! &\phi_1'' (t_{n-1}')\\
&\vdots & {\vdots} &\vdots\\
&\phi_{n-1}''(t_1) &\phi_{n-1}''(t_2') \ \   \cdots &\phi_{n-1}''(t_{n-1}')\end{matrix}\right|
$$
\medskip

and the non-vanishing can again be ensured by \eqref{2.1}.

\section
{Decoupling inequality for curves}

Next, we use (1.3) to derive a decoupling inequality for curves (a variant of this approach appears in \cite{B-D2}.
 
Let $\Gamma_1, \ldots, \Gamma_{n-1} \subset \Gamma\subset \mathbb R^n$ be as in \S2.
Let $\delta>0$ and denote by $\Gamma_j^\delta$ a $\delta$-neighborhood of $\Gamma_j$.
\bigskip

Assume $\supp\widehat{f_j}\subset\Gamma_j^\delta$.

\bigskip
Write with $x=(x_0, x_1, \ldots, x_{n-1})\in \mathbb R^n$ and $\Phi$ as above
\begin{align}\label{3.1}
&\prod^{n-1}_{j=1} \Big[\int_{I_j} \widehat {f_j} (t_j) e (\langle\Phi (t_j). x\rangle) dt_j\Big]= \nonumber\\
&\int_{I_1} \cdots\int_{I_{n-1}}\Big[\prod^{n-1}_{j=1} \widehat {f_j}(t_j)\Big]
e\big(\langle\Phi(t_1)+\cdots+\Phi(t_{n-1}). x\rangle\big) dt_1 \ldots
dt_{n-1}=\nonumber\\
&\int_S \Big[\prod^{n-1}_{j=1} \widehat{f_j}(t_j)\Big] e(\xi. x) \Omega(\xi)d\xi
\end{align}
with $\Omega$ some smooth density on $S$.

Let $p=\frac{2(n+1)}{n-1}$ and apply the decoupling inequality for $S$ stated in \eqref{1.3} of Section 1.
Observe that by the regularity of $D_1$ in \eqref{2.5}, a partition of $S$ in $\sqrt\delta$-caps $\tau_\alpha \subset S$ is equivalent to a
partition of $I_1\times \cdots \times I_{n-1}$ in $\sqrt\delta$-cubes.
Hence, denoting by $J\subset [0, 1]$ $\sqrt\delta$-intervals, we obtain
\begin{align}\label{3.2}
&\Vert\eqref{3.1}\Vert_{L^p(B_{\frac 1\delta})}\ll\nonumber\\
&\delta^{-\frac{n-1}{2(n+1)}-\ve} 
\Big\{\sum_{\substack
 {J_1, \ldots, J_{n-1}\\ J_j\subset I_j}}
 \Big\Vert\prod^{n-1}_{j=1} \Big[\int_{I_j} \widehat{f_j} (t_j) e\big(t_jx_0+f_1(t_j)x_1
+\cdots+ f_{n-1}(t_j) x_{n-1}\big) dt_j\Big] \Big\Vert^p_{L^p(B_{\frac 1\delta})}\Big
\}^{1/p}.
\end{align}

Next take $N=\frac 1\delta$ and discretize inequality \eqref{3.2}  by setting $t=\frac kN, k\in \big\{\frac N2, \ldots, N\}$.

\bigskip
This leads to the following inequality for separated intervals $U_1, \ldots, U_{n-1}\subset\{\frac N2, \ldots, N\}$
\begin{align}\label{3.3}
&\Big\Vert\prod^{n-1}_{j=1} \Big|\sum_{k\in U_j} a_k e \Big(kx_0+N\vp_1\Big(\frac kN\Big) x_1+\cdots + N\vp_{n-1} \Big(\frac kN\Big) x_{n-1}\Big)
\Big|\Big\Vert_{L^p([0, 1]^n)}\ll\nonumber\\
& N^{\frac {n-1}{2(n+1)}+\ve}\Big(\sum_{\substack{V_1, \ldots, V_{n-1}\\ V_j\subset U_j}}\Big\Vert\prod^{n-1}_{j=1} \Big|\sum_{k\in V_j} a_k
e(\cdots) \Big| \, \Big\Vert^p_{L^p([0, 1]^n)}\Big)^{\frac 1p}
\end{align}
with $V\subset\{\frac N2, \ldots, N\}$ running in a partition in $\sqrt N$-size intervals.
\bigskip

Note that the domain $[0, 1]^n$ may always be replaced by a larger box $\prod^{n-1}_{j=0}[0, K_j], K_j\geq 1$.
In particular, the function
$$
kx_0 +N\phi_1\Big(\frac kN\Big) x_1+\cdots+N\phi_{n-1} \Big(\frac kN\Big) x_{n-1}
$$
in \eqref{3.3} may be replaced by
$$
kx_0+N_1\phi_1\Big(\frac kN\Big) x_1+\cdots+ N_{j-1} \phi_{n-1} \Big(\frac kN\Big) x_{n-1}\text { where } \  N_1, \ldots, N_{j-1}\geq N.
$$

Take $\phi_1(t) =t^2, N_1= N^2, N_2 =\cdots= N_{j-1}=N$.
We obtain 
\begin{align}\label{3.4}
&\Big\Vert \prod^{n-1}_{j=1} \Big|\sum_{k\in U_j} e(kx_0+k^2 x_1+N\phi_2\Big(\frac kN\Big) x_2+\cdots +N\phi_{n-1} \Big(\frac kN\Big)
x_{n-1}\Big)\Big|\Big\Vert_{L^p([0, 1])^n)}\nonumber\\
&\ll N^{\frac {n-1}{2(n+1)} +\ve} \Big(\sum_{V_1, \ldots, V_{n-1}} \Big\Vert\prod^{n-1}_{j=1} \Big|\sum_{k\in V_j}
e(\cdots)\Big|\Big\Vert^p_{L^p([0, 1]^n)} \Big)^{\frac 1p}.
\end{align}
\bigskip

Our next task is to bound the individual summands in \eqref{3.4}.
\bigskip

Write $\bar k=(k_1, \ldots, k_{n-1})\in V_1 \times \cdots \times V_{n-1}$ as $\bar k=\bar\ell+\bar m$, $\ell_j$ the center of $V_j$ and
$|m_j|<\sqrt N$. Hence
\begin{align}\label{3.5}
&\sum^{n-1}_{j=1} \Big(k_jx_0+k_j^2 x_1+N\phi_2 \Big(\frac {k_j}N\Big) x_2+\cdots+ N\phi_{n-1} 
\Big(\frac {k_j}N\Big) x_{n-1}\Big)=\nonumber\\
&m_1 \Big(x_0+2\ell_1 x_1 +\phi_2'\Big(\frac {\ell_1}N\Big) x_2 +\cdots+ \phi_{n-1}'\Big(\frac {\ell_1}N\Big)x_{n-1}\Big)+\nonumber\\
&\vdots\nonumber\\
&+ m_{n-1} \Big(x_0+2\ell_{n-1} x_1+ \phi_2' \Big(\frac {\ell_{n-1}}N\Big) x_2+\cdots+ \phi_{n-1}' \Big( \frac{\ell_{n-1}} N\Big)
x_{n-1}\Big)+\nonumber\\
&(m_1^2 +\cdots + m^2_{n-1})x_1 +\psi(\bar m, x)
\end{align}
where $|\psi(\bar m, x)|<o(1)$ and $|\partial_m\psi(\bar m, x)|<O(N^{-\frac 12})$ since $|\bar m|<\sqrt N$ and $|x|<1$.

Thus $\psi(\bar m, x)$ may be dismissed in \eqref{3.4} when evaluating \hfill\break
\bigskip
$\Vert\prod^{n-1}_{j=1} \big|\sum_{k\in V_j}e(\cdots)
\big|\big\Vert_{L^p([0, 1]^n)}$. Make an affine change of variables
$$
\begin{pmatrix} y_1\\ \vdots\\ y_{n-1}\end{pmatrix} 
=A\begin{pmatrix} x_0\\ x_2\\ \vdots\\ x_{n-1}\end{pmatrix} \text { with }
A= \left[\begin{pmatrix} &1 &\phi_2' (\frac {\ell_1}N)&\cdots&\phi_{n-1}'(\frac{\ell_1}N)\\
&\vdots &\vdots && \vdots&\\
&1&\phi_2'(\frac{\ell_{n-1}}N)&\cdots& \phi_{n-1}' (\frac{\ell_{n-1}}N)\\
\end{pmatrix}\right]
$$
in the $(x_0, x_2, \ldots, x_{n-1})$ variables, noting that this linear coordinate change can be assumed regular provided
$W(\phi_2'', \ldots, \phi_{n-1}'')\not= 0$ (which is implied by \eqref{2.1} for $\phi_1(t)=t^2)$.

Next, using periodicity, another coordinate shift leads to
\begin{align}\label{3.6}
&\big\Vert\prod^{n-1}_{j=1} \Big|\sum_{k\in V_j} e(\cdots) \Big|\Big\Vert_{L^p([0, 1]^n)}\sim\nonumber\\
&\Big\Vert\sum_{m_1, \ldots, m_{n-1} <\sqrt N} 
e\big(m_1y_1+\cdots+ m_{n-1} y_{n-1}  +
(m^2_1+\cdots+ m^2_{n-1})x_1\big) \Big\Vert_{L^p_{x_1, y_1, \ldots, y_{n-1}}([0, 1]^n)}\nonumber\\
&\lll N^{\frac {n-1}4+\ve}
\end{align}
by the Strichartz inequality on $\mathbb T^n$.

Summarizing, we proved the following multi-linear mean value theorem.

\begin{theorem}\label{Theorem1}
Assume $n\geq 3$ and $\vp_2, \ldots, \vp_{n-1}$ satisfying
$$
W(\vp_2''', \ldots, \vp_{n-1}''')\not= 0.
$$
Let $U_1, \ldots, U_{n-1} \subset \big[\frac N2, N]\cap\mathbb Z$ be $O(N)$-separated intervals. Then
\begin{align}\label{3.7}
&\Big\Vert\prod^{n-1}_{j=1} \Big|\sum_{k\in U_j} e\Big(kx_0+k^2 x_1+N\vp_2\Big(\frac kN\Big)x_2+\cdots+ N\vp_{n-1} \Big(\frac
kN\Big) x_{n-1}\Big)\Big|\Big\Vert_{L^{\frac {2(n+1)}{n-1}}([0, 1]^n)}\leq\nonumber\\
&N^{\frac {n-1}2+\ve}.
\end{align}
\end{theorem}

\noindent
{\bf Remarks.}

(i) Theorem 1 remains valid (following the same argument) with coefficients $a_k, k\in U_j$ and r.h.s. replaced by
$\prod^{n-1}_{j=1} (\sum_{k\in U_j}|a_k|^2)^{\frac 12}$.

(ii) Note that \eqref{3.7} is best possible.  Indeed, restricting $|x_0|< \frac 1N, |x_1|<\frac 1{N^2}$, $|x_2|<\frac 1N,
\ldots, |x_{n-1}|<\frac 1N$, one gets the contribution
$$
N^{n-1-\frac{n-1}{2(n+1)} (n+1)} =N^{\frac {n-1}2}.
$$

(iii) Also, as we will see shortly, \eqref {3.7} is only valid in the above multi-linear form.
\bigskip

\section
{Mean values estimates for the 8th and 10th  moment}

Note that \eqref{3.7} is the multi-linear version of an estimate on
\be\label{4.1}
\Big\Vert \sum_{k\sim N} e\Big(k x_0+k^2x_1+N\vp_2 \Big(\frac kN\Big) x_2+      \cdots+ N\vp_{n-1} \Big(\frac kN\Big)
x_{n-1}\Big)
\Big\Vert_{L^{2(n+1)}([0, 1]^n)}.
\ee

Denoting $f_I=\sum_{k\in I} e\big(k x_0+k^2x_1+ N\vp_2\big(\frac kN\big) x_2+   \cdots+ N\vp_{n-1} \big(\frac kN\big)
x_{n-1}\big)$ for $I\subset [1, N]$  a subinterval, one adopts the
following argument from \cite{B-G}.
Partition $[1, N]$ in intervals $I$ of size $N^{1-\tau}$ ($\tau>0$ small).
Fix a point $x$ and distinguish the following  two scenarios. 
Either we can find $n-1$ intervals $I_1, \ldots, I_{n-1}$ that are $O(N^{1-\tau})$-separated and such that
\be\label{4.2}
|f_{I_j} (x)|>N^{-2\tau} |f(x)| \ \text { for } \ 1\leq j\leq n-1
\ee
or for some interval $I$, we have
\be\label {4.3}
|f_I (x)|>c|f(x)|.
\ee

The contribution of \eqref{4.2} is captured by the multi-linear estimate \eqref{3.7} and we obtain $N^{\frac 12+c\tau}$.
For the \eqref{4.3}-contribution, bound by
$$
\max |f_I|\leq \Big[\sum_{I} |f_I|^{2(n+1)}\Big]^{\frac 1{2(n+1)}}
$$
contributing to
\be\label{4.4}
\Big[\sum_I \Vert f_I\Vert^{2(n+1)}_{2(n+1)}\Big]^{\frac 1{2(n+1)}}.
\ee
One may then repeat the process to each $f_I$.
Note that after a coordinate change in $x_0, x_1$, we obtain exponential sums of the form
\be\label{4.5}
F(x)=\sum_{\ell \sim  M} e\Big(\ell x_0+\ell^2x_1+N\vp_{_2} \Big(\frac kN+\frac\ell N\Big) x_2+\cdots+ N\vp_{n-1}
\Big(\frac kN+\frac \ell N\Big) x_{n-1}\Big)
\ee
with $M=N^{1-\tau}$, $k\sim N$ fixed.
Set for $j=1, \ldots, n-1$.
\be\label{4.6}
N\vp_j\Big(\frac kN+\frac \ell N\Big)= M\psi_j\Big(\frac\ell M\Big)
\ee
with
\be\label{4.7}
\psi_j(t)=\frac NM \vp\Big(\frac kN+\frac MNt\Big).
\ee
However the Wronskian condition $W(\psi_2''', \ldots, \psi_{n-1}''')> O(1)$ deteriorated.
For $n=3$ we will nevertheless be able to retrieve easily the expected bound, while for $n\geq 4$, the linear
bounds turn out to be weaker than the multi-linear one.
\bigskip

Let $n=3$.
Then $\psi_2'''(t) =\frac {M^2}{N^2} \vp'''(\frac kN+\frac MN t)=O(\frac {M^2}{N^2})$
and replacing $\psi_2= \frac{M^2}{N^2} \tilde\vp_2, x_2' =\frac {M^2}{N^2} x_2$, this leads to
$$
\begin{aligned}
\Vert F\Vert_{L^8_{x_0, x_1, x_2=O(1)}}&
\sim \Big(\frac NM\Big)^{\frac 14}\Big\Vert\sum_{\ell\sim M} e\Big(\ell x_0+\ell^2 x_1
+M\tilde \vp_2 \Big(\frac\ell M\Big) x_2'\Big)
\Big\Vert_{\substack{{L^8 _{x_0, x_1=O(1)}}\\  {x_2'=O(\frac {M^2}{N^2})}}}\\
&\leq \Big(\frac NM\Big)^{\frac 14} \Vert \cdots\Vert_{L^8_{x_0, x_1, x_2'=O(1)}}\\
&<\Big(\frac NM\Big)^{\frac 14} M^{\frac 12+\ve}
\end{aligned}
$$
assuming the expected bound at scale $M$.
The bound on \eqref{4.4} becomes then
$$
\Big(\frac NM\Big)^{\frac 38} M^{\frac 12+\ve} =\Big(\frac MN\Big)^{\frac 18} N^{\frac 12+\ve}
$$
and summing over dyadic $M<N$ we reproved the main result from \cite{B-I2}. 

\begin{theorem}\label{Theorem2}
{\rm \cite{B-I2}}.

Assume $\vp''' \not= 0$. Then
\be\label{4.8}
\Big\Vert\sum_{k\sim N} e\Big( kx_0+ k^2x_1+ N\vp \Big(\frac kN\Big) x_2\Big)\Big\Vert_8 \ll N^{\frac 12+\ve}.
\ee
\end{theorem}
\bigskip

Note that in their application to $\zeta(\frac 12+it)$, $\vp(t)=t^{3/2}$.
\bigskip

It is interesting to note that unlike in \cite{B-I2}, our derivation of \eqref{4.8} did not make use of Poisson summation
(i.e. Process B).

The work of \cite{B-I1} was extensively refined by Huxley and his
collaborators, resulting in his book \cite{H}.

The present discussion is relevant to the so called `First Spacing Problem'; \eqref {4.8} indeed means that the 
system
$$
(4.9)
\begin{cases}
k_1+k_2+k_3+k_4 = k_5+\cdots+k_8\\
k_1^2+\cdots +k_4^2= k^2+\cdots +k^2_8\\
k_1^{3/2} +\cdots+ k_4^{3/2} =k_5^{3/2} +\cdots+k_8^{3/2}+ O(\sqrt N).
\end{cases}
$$
has at most $N^{4+\ve}$ solutions in integers $k_1, \ldots, k_8\sim N$
(the statement is clearly optimal).
\bigskip

Huxley considers the more elaborate problem in 10-variables
$$
(4.10) 
\begin{cases}
k_1+k_5=k_6+\cdots+k_{10}\\
k_1^2+\ldots+k^2_5= k_6^2+\cdots+k^2_{10}\\
k_1^{3/2}+\cdots+ k_5^{3/2}= k_6^{3/2}+\cdots+ k_{10}^{3/2}+O(\delta N^{3/2})\\
k_1^{1/2}+\cdots+ k_5^{1/2}= k_6^{1/2}+\cdots + k_{10}^{1/2}+O(\Delta N^{\frac 12})
\end{cases}
$$
(see \cite{H}, \S11) for which the number $N_{10}(\delta, \Delta)$ of solutions is given by the 10th moment  
\stepcounter{equation}
\stepcounter{equation}
\be\label{4.11}
\Big\Vert\sum_{k\sim N} e\Big(kx_0+k^2x_1+\frac 1\delta \Big(\frac kN\Big)^{3/2} x_2+\frac 1\Delta
\Big(\frac kN\Big)^{1/2} x_3\Big)\Big\Vert_{L^{10}_{x_0, x_1, x_2, x_3}}^{10}.
\ee

In the applications to exponential sums, $\Delta =\delta N$, $\frac 1{N^2} <\delta <\frac 1N$.
In this setting, the following key inequality appears in \cite{H-K}.

\begin{theorem}\label{Theorem3} \cite{H-K}. With $\Delta =\delta N, \frac 1{N^2} <\delta<\frac 1N$, we have
\be \label{4.12}
N_{10}(\delta, \delta N) \ll \delta. N^{7+\ve}.
\ee
\end{theorem}
\medskip

In what follows, we estimate \eqref{4.11} using Theorem 7 and will in particular retrieve \eqref{4.12} in a stronger form.

\medskip
Start by observing that, as a consequence of \eqref{3.7}, for $U_1, U_2, U_3$ and $\vp_2, \vp_3$ as in Theorem~1
\begin{align}\label{4.13}
\int^1_0\int^1_0\int_0^1\int^1_0 &\Big\{\prod^3_{j=1} \Big|\sum_{k\in U_j} e(kx_0+k^2x_1+\frac 1\delta
\vp_2\Big(\frac kN\Big)x_2+ \frac 1\Delta \vp_3\Big(\frac kN\Big) x_3\Big|^{\frac {10}3}\Big\} dx_0dx_1 dx_2 dx_3\ll\nonumber\\
& [\min(\delta N, N)+1] [\min(\Delta N, N)+1]N^{5+\ve}
\end{align}

Using the scale reduction described in \eqref{4.1}-\eqref{4.7}, we also need to evaluate the contributions of
\stepcounter{equation}
\be\label{4.15}
\frac NM\cdot (4.14)
\ee
with
$$
(4.14) =\int^1_0\int^1_0\int^1_0\int^1_0\Big\{\prod^3_{j=1} \Big|\sum_{\ell \in U_j} e(\ell x_0+\ell^2 x_1+\frac 1\delta 
\vp_2\Big(\frac {k+\ell}N\Big)x_2
+\frac 1\Delta \vp_3\Big(\frac {k+\ell}N\Big) x_3\Big|^{\frac {10}3}\Big\}
$$
where $k\in [\frac N 2, N]$, $I=[k, k+M[$ and $U_1, U_2, U_3$ are $\sim M$ separated subintervals of size $\sim M$ in $I$.
By a change of variables in $x$, the phase function in \eqref{4.16} may be replaced by
\be\label{4.16}
e\Big(\ell x_0+\ell^2 x_1+\frac {M^3}{\delta N^3} \tilde\vp_2 \Big(\frac \ell M\Big) x_2+\frac {M^4}{\Delta N^4}\tilde\vp_3 \Big(\frac \ell M\Big) x_3\Big)
\ee
where $\tilde\vp_2(t)$ has leading monomial $t^3$ and $\tilde\vp_3(t)$ leading monomial $t^4$.
Hence $W(\tilde\vp_2''', \tilde\vp_3''')>c$ and \eqref{4.13} is applicable to (4.14) with $N, \delta, \Delta$ replaced by $M, \frac {\delta N^3}{M^3}, \frac
{\Delta N^4}{M^4}$.
Therefore
\be\label{4.17}
\eqref{4.15} \ll \Big[1+\min\Big(\frac {\delta N^3}{M^2}, M\Big)\Big] \Big[1+\min \Big(\frac {\Delta N^4}{M^3}, M\Big)\Big] M^4 N^{1+\ve}
\ee
and \eqref{4.17} needs to be summed over dyadic $M<N$.  One easily checks that the conclusion is as follows

\begin{theorem}\label {Theorem4}

Assume $W(\vp_2''', \vp_3''')\not= 0$ and $\delta<\Delta$. Then
$$
\begin{aligned}
&\Big\Vert \sum_{k\sim N} e\Big(kx_0+k^2x_1+\frac 1\delta \vp_2\Big(\frac kN\Big)x_2 +\frac 1\Delta \vp_3 \Big(\frac kN\Big)
x_3\Big)\Big\Vert_{10}^{10} \ll\\
&[\delta\Delta^{3/4} N^7 +(\delta+\Delta)N^6+N^5]N^\ve.
\end{aligned}
\eqno{(4.18)}
$$
In particular
$$
N_{10} (\delta, \Delta)< (4.18).
$$
\end{theorem}

Hence, we are retrieving Theorem 3.
\medskip

\noindent
{\bf Remark.}
We make the following comment on the role of the first term in the r.h.s. of (4.18), relevant to the Remark following
Theorem 1.

Partition $[\frac N2, N]$ in intervals $I=[n, n+M]$ of size $M$.
Obviously $N_{10}(\delta, \Delta)$ is at least $\frac NM$ times a lower bound on the number of solutions of
$$
\begin{cases}
\left. \begin{matrix} m_1+\cdots+ m_5=m_6+\cdots+m_{10}\\
m_1^2+\cdots+ m^2_5=m^2_6+\cdots+m_{10}^2\end{matrix}
\quad\qquad\qquad \, \qquad\qquad\qquad\qquad\qquad \right\}{(4.19)}\\
\left.\begin{matrix} \Big(\frac {n+m_1} N\Big)^{\frac 32} +\cdots+ \Big(\frac {n+m_5}{N}\Big)^{\frac 32} =
\Big(\frac {n+m_6}N\Big)^{\frac 32}+\cdots +\Big(\frac {n+m_{10}} N\Big)^{\frac 32}+O(\delta)\\
\Big(\frac{n+m_1}N\Big)^{\frac 12} +\cdots+ \Big(\frac{n+m_5}N\Big)^{\frac 12}=
\Big(\frac {n+m_6}N\Big)^{\frac 12}+\cdots+ \Big(\frac {n+m_{10}}N\Big)^{\frac 12}+O(\Delta)\end{matrix}\right\}
{(4.20)}
\end{cases}
$$
Since
$$
\begin{aligned}
&\Big(\frac {n+m}N\Big)^{\frac 32} =\Big(\frac nN\Big)^{\frac 32}+\frac 32\Big(\frac nN\Big)^{\frac 12} \frac mN+\frac
32\Big(\frac nN\Big)^{-\frac 12}\Big(\frac mN\Big)^2 -\frac 1{16} \Big(\frac nN\Big)^{-\frac 32} \Big(\frac mN\big)^3+\ldots\\
&\Big(\frac{n+m}N\Big)^{\frac 12}=\Big(\frac nN\Big)^{\frac 12} +\frac 12\Big(\frac nN\Big)^{-\frac 12} \frac mN -\frac 18\Big(\frac
nN\Big)^{-\frac 32} \Big(\frac mN\Big)^2+\frac 3{16} \Big(\frac nN\Big)^{-\frac 52}\Big(\frac mN\Big)^3 -\frac {15}{128} \Big(\frac
nN\Big)^{-\frac 72}
\Big(\frac mN\Big)^4+\ldots
\end{aligned}
$$
the equations (4.20) may be replaced by
$$
\begin{cases}
\vp\Big(\frac {m_1}N\Big)+\cdots - \vp\Big(\frac{m_{10}}N\Big)=O(\delta)\\
\psi \Big(\frac {m_1}N\Big)+\cdots - \psi\Big(\frac {m_{10}}N\Big) =O(\Delta)
\end{cases}
$$
with $\vp, \psi$ of the form $\vp(t)=a_3t^3+a_4 t^4+\cdots$ and $ \psi (t)= b_3 t^3 + b_4 t^4+\cdots$ and where $\left|\begin{matrix}
a_3&b_3\\ a_4&b_4\end{matrix}\right|\not= 0$.

\stepcounter{equation}
\stepcounter{equation}

Assume $\delta < \Delta$ and replace $\psi$ by $\psi_1=\psi -\frac {b_3}{a_3}\vp = c_4 t^4 +\cdots$
Writing $\frac mN=\frac MN\frac mM$, we obtain conditions of the form
\stepcounter{equation}
\be\label{4.21}
\begin{cases} \tilde\vp\Big(\frac {m_1}M\Big)+\cdots -\tilde\vp\Big(\frac {m_{10}}N\big)< O\Big(\frac {N^3}{M^3}\delta\Big)\\
\tilde\psi_1\Big(\frac {m_1}M\Big)+\cdots-\tilde\psi_1\Big(\frac {m_{10}}N\Big)<O\Big(\frac {N^4}{M^4}\Delta\Big)\end{cases}
\ee
where $\tilde\vp = t^3+\cdots, \tilde\psi_1 =t^4+\cdots$.
Consider the system (4.19)+(4.21) with $m_i\leq M$.
Clearly the number of solutions is at least
$$
M^7 \min\Big(1, \frac {N^3}{M^3}\delta\Big) .\min\Big( 1, \frac {N^4}{M^4}\Delta\Big).
$$
Taking $M= \Delta^{1/4}N$,  we obtain $N^7\delta\Delta$.
The quantity is multiplied further with $\frac NM$, leading to a lower bound 
$\delta\Delta^{3/4}N^7$ for
$N_{10}(\delta, \Delta)$.

This shows that the first term in (4.18) (apart from the $N^\ve$ factor) is also a lower bound.
\bigskip

In our applications, $\Delta$ tends to be much larger then $\delta$ which makes $\Delta N$ the leading term in (4.18).
Next, we develop an argument to reduce the weight of $\Delta N$ by
involving also some ideas and techniques from \cite{H}.
It is likely that our presentation can be improved at this point.

We will need the following variant of van der Corput's exponential sum bound (cf. [Ko], Theorem 2.6).

\begin{lemma}\label{Lemma5}
Assume $f$ a smooth function on $I=[\frac N2, N]$ and $f^{(3)} \sim\lambda_3$.
Let $\{ V_j\}$ denote a partition of $I$ in intervals of size $D$. Then
$$
\sum_j\Big|\sum_{n\in V_j} e\big(f(n)\big)\Big|^2 \lesssim
\left\{
\begin{aligned} &N+D^{\frac 12} \lambda_3^{-\frac 12} +D^{\frac 32} \lambda_3^{\frac 12} N\qquad\qquad\qquad\qquad& (4.22)\\ 
&ND\lambda_3^{\frac 13}+D\lambda_3^{-\frac 13} \text { if } \ D>\lambda_3^{-\frac 13}.& (4.23) 
\end{aligned}\right.
$$
\end{lemma}

\stepcounter{equation}

We first proceed with a multi-linear bound
considering instead of \eqref{4.13} 5-linear expressions with $U_j\subset [{\frac N2, N}]$ of size $\sim N$ and $\sim N$
separated $(1\leq j\leq 5)$
\stepcounter{equation}
\be\label{4.24}
\int\Big\{ \prod^5_{j=1} \Big|\sum_{k\in U_j} e\Big(kx_0+k^2x_1+\frac 1\delta\Big(\frac kN\Big)^{3/2} x_2+\frac 1\Delta \Big(\frac
kN\Big)^{\frac 12} x_3\Big)\Big|^2\Big\} dx_0dx_1 dx_2 dx_3.
\ee
This quantity will increase by increasing $\delta$ and we replace $\delta$ by a parameter $\delta_1>\delta$ to be specified.
An application of H\"older's inequality permits then to bound \eqref{4.24} by \eqref{4.13} with $\delta$ replaced by $\delta_1$.

Assuming  $\Delta<1$, perform a decoupling at scale $N\Delta^{\frac 12}$ using \eqref{3.2}.
This gives an estimate on the l.h.s. of \eqref{4.13} by
\be\label{4.25}
\Delta^{-1}\int \prod^3_{j=1} \Big[\sum_{V_j\subset U_j} \Big|\sum_{k\in V_j} e(\cdots) \Big|^{\frac {10} 3}\Big] dx
\ee
with $V_j\subset U_j$ a partition in $N\Delta^{\frac 12}$-intervals.
Using again H\"older's inequality, one may bound
\be\label{4.26}
\prod^3_{j=1}\Big[\sum_{V_j\subset U_j}\Big| \quad \Big|^{\frac {10}3}\Big] \leq \Big(\sum_{V\subset [\frac N2, N]} \Big| \quad \Big|^2\Big)
\prod^2_{j=1}  \Big(\sum_{V_j\subset U_j}\Big| \quad\Big|^4 \Big)  + \cdots
\ee
where $\cdots$ refers to the pairs $U_2, U_3$ and $U_3, U_1$ instead of $U_1, U_2$.

Specifying in \eqref{4.13}, with $\delta$ replaced by $\delta_1$, a range
\be\label{4.27}
x_2\sim X_2 <1 \text { assuming } \ X_2\Delta>100\delta_1
\ee
an application of (4.22) to the first factor of \eqref{4.26} with $D=\Delta^{\frac 12} N, \lambda_3\sim \frac {X_2}{\delta_1 N^3}$ 
gives the bound
\be\label{4.28}
N+\delta_1^{\frac 12} \Delta^{\frac 14} N^2 X_2^{-\frac 12}+\delta_1^{-\frac 12} \Delta^{\frac 34}N X_2^{\frac 12}.
\ee
We always assume
\be\label{4.29}
\Delta N>100
\ee
(this condition remains clearly preserved at lower scales, cf. \eqref{4.17}).

Apply the bilinear estimate (Theorem 1 with $n=3$) to the second factor of (4.26) considering the variables $x_0, x_1, x_2$ and
restricting $x_2\sim X_2$.
By \eqref{4.24}, this gives the contribution
$$
(1+\delta_1^{\frac 12} \Delta^{\frac 14} NX_2^{-\frac 12}+
\delta_1^{-\frac 12}\Delta^{\frac 34} X_1^{\frac 12})
(X_2+\delta_1 N)N^{5+\ve}.
$$
Assuming $X_2>\delta_1N$, which by \eqref{4.29} implies \eqref{4.27}, gives the bound
$$
(1+\delta_1^{\frac 12} \Delta^{\frac 14} N+\delta_1^{-\frac 12}\Delta^{\frac 34})N^{5+\ve}.
$$
The contribution of $X_2<\delta_1N$ is estimated by
$$
\begin{aligned}
\Delta^{-1}&\int^1_0\int_0^1\int_0^{\delta_1N}\int_0^1
\Big\{\prod^3_{j=1} \Big[\sum_{V_j\subset U_j}\Big| \quad\Big| ^{\frac {10} 3} \Big]\Big\} dx_0 dx_1 dx_2 dx_3\\
&\ll \Delta^{-1} \Big(\frac N{\Delta^{1/2}N}\Big)^3(\Delta N)\delta_1 N (\Delta^{\frac 12}N)^5 \ll \Delta\delta_1 N^{7+\ve}
\end{aligned}
$$
assuming $\delta_1N<1$.
This gives
\be\label{4.30}
(1+\delta_1^{-\frac 12}\Delta^{\frac 34} +\delta_1^{\frac 12} \Delta^{\frac 14} N+\Delta\delta_1N^2) N^{5+\ve}
\ee
and setting $\delta_1=\Delta^{\frac 12} N^{-1}$, assuming $\delta_1>\delta$ gives
$$
\big((\Delta N)^{\frac 12}+ \Delta^{3/2}N\big) N^{5+\ve}.
$$
If $\delta_1\leq  \delta$, use \eqref{4.30} with $\delta_1=\delta$.

Thus the multi-linear contribution in the 10th moment may be estimated by
\be\label{4.31}
\Big(\delta^{\frac 12}\Delta^{\frac 14}N+\Delta \delta N^2+(\Delta N)^{\frac 12}+\Delta^{3/2} N\Big)
N^{5+\ve}.
\ee

Next, consider te lower scale contributions
\be\label{4.32}
\sum_{\substack{I\subset [\frac N2, N]\\ {|I|=M}}}\Big\Vert\sum_{n\in  I} e(\cdots)\Big\Vert^{10}_{10}.
\ee

Fixing $M<N$ and replacing $\delta$, resp. $\Delta$, by $\frac {N^3}{M^3}\delta$, resp. $\frac {N^4}{M^4}\Delta$, we obtain the
bound
\be\label{4.33}
\big(\delta^{\frac 12}\Delta^{\frac 14} N+(\Delta N)^{\frac 12}\big) N^{5+\ve} +(\Delta\delta N+\Delta^{3/2})\frac {N^{7+\ve}}M
\ee
for the multi-linear contribution at scale $M$.

On the other hand, we can also make a crude estimate using the $L^8$-norm, leading to the contribution
\be\label{4.34}
\frac NM M^2\Big(1+\delta\frac {N^3}{M^2}\Big) M^{4+\ve}\ll NM^{5+\ve}+\delta N^4 M^{3+\ve}
\ee
and $(4.34) \ll (\Delta N)^{\frac 12} N^{5+\ve}$ for $M<(\Delta N)^{\frac 1{10}} N^{4/5}$ and $\delta<N^{-7/5}(\Delta N)^{\frac
15}$.

Hence we get

\begin{lemma}\label{Lemma6}
For $N\Delta>1$ and $\delta<N^{-7/5}(\Delta N)^{\frac 15}$
\be\label{4.35}
N_{10} (\delta, \Delta)\ll \big(\delta^{\frac 12} \Delta^{\frac 14} N+(\Delta N)^{\frac 12}\big)N^{5+\ve}
+ (\Delta\delta N+\Delta^{\frac 32})(\Delta N)^{-\frac 1{10}} N^{\frac {31}5+\ve}.
\ee
\end{lemma}

Next, recall Lemma 11, 3.3 in \cite{H}.

\begin{lemma}\label{Lemma7}
Assume $\frac 1N>\delta>\frac 1{N^2}$ and $\frac 1N<\Delta<\delta N$.
Letting
$$
2\leq T\leq \frac 1{\sqrt {\delta N}}
$$
be a parameter, the following inequality holds
\be\label{4.36}
N_{10} (\delta,\Delta)\lesssim \frac 1T N_{10} (T^2\delta, T\Delta)+N_{10}(\delta, CT\delta).
\ee
\end{lemma}

Combining Theorem 4, Lemmas 6 and 7 (applied with $T=\Delta N)$ gives

\begin{lemma}\label{Lemma8}
Assume $\frac 1N>\delta>\frac 1{N^2}, \frac 1N<\Delta<\delta N$ and $\Delta\sqrt\delta N^{3/2}<1$.
Then
\be\label{4.37}
N_{10}(\delta, \Delta)<\big(1+\delta (\Delta N)^{\frac 85} N^{\frac 75}\big)
\big(1+\delta^{\frac 12}(\Delta N)^{\frac 12} N^{\frac 34})N^{5+\ve}+ \big(\delta(\Delta N)^{\frac {14}5}+
(\Delta N)^{\frac 3{10}}\Delta^{\frac 32}\big)N^{\frac {31}5+\ve}.
\ee
\end{lemma}

Setting $\Delta=\delta N$ leads to the following strengthening of Theorem 3

\begin{theorem}\label{Theorem9}
For $N^{-\frac{33}{18}}\geq\delta\geq N^{-2}$, we have $N_{10} (\delta, \delta N)\ll N^{5+\ve}$.
\end{theorem} 

Note that in view of the Remark following Theorem 4, the conclusion of Theorem 9 fails for $\delta>N^{-11/7}$.

\section
{On an inequality of Robert and Sargos}

In \cite{R-S} established the inequity
\be\label{5.1}
I_{10} =\int^1_0\int^1_0 \Big(\sum_{n\sim N} e(n^2 x+n^4y)|^{10} dxdy \ll N^{\frac {49}8+\ve}
\ee
which they applied to obtain new estimates on Weyl sums.
An improvement of \eqref{5.1} appears in \cite{P}, who obtains
\be\label{5.2}
I_{10} \ll N^{6+\ve}.
\ee
Using our methods, we present a further improvement.

\begin{theorem}\label{Theorem11}
\be\label{5.3}
I_{10} \ll N^{\frac {17}3+\ve}
\ee
\end{theorem}

The corresponding improvement in Weyl's inequality following Heath-Brown's method was recorded in the Introduction.

Note that bounding $I_{10}$ is tantamount to estimating the number of integral solutions $n_i\sim N$ $(1\leq 1\leq 10)$ of the
system
\be\label{5.4}
\begin {cases}
n_1^2+n^2_2+n^2_3+n_4^2 +n_5^2 = n_6^2+n_7^2+n_8^2+n_9^2+n_{10}^2\\
n_1^4+n_2^4+n_3^4+n_4^4+n_5^4 =n_6^4+n_7^4+n_8^4 +n_9^4+n_{10}^4.
\end{cases}
\ee

The problem is not shift invariant and therefore as it stands not captured by a Vinogradov mean value theorem of the usual kind.
Following Wooley's approach for $(n, n^3)$ (see \cite{W}), knowledge of the (conjectural) optimal VMVT for $k=4$ (which would involve
the 20th moment) and interpolation with the 6th moment would at the best deliver
$I_{10}\ll N^{\frac {41}{7}}$, inferior to \eqref{5.3}.

A crude summary of our argument.
As in \cite{R-S}, we need to consider the more general expressions
\be\label{5.5}
I_p(\lambda)=\int^1_0\int^1_0 \Big|\sum_{n\sim N} e(n^2 x+\lambda n^4 y)\Big|^p dxdy
\ee
with $p\geq 6$ and $0<\lambda\leq 1$.
A first step is an application of the decoupling theorem from \cite{B-D} for planar
curves similarly as in \cite{B-D}, Theorem 2.18 (where an
extension of the result $I_6(N^{-3})\ll N^{3+\ve}$ from \cite{R-S} is established).
At this stage, one gets shorter sums, of length $M$ say, i.e. $n\in [n_0, n_0+M]$ with $n_0$ ranging in $\big[N, \frac N2\big]$.
Exploiting $n_0$ as an additional variable leads then to mean value expressions of the form
\be\label{5.6}
\int_0^1\int_0^1\int_0^1\int^1_0\Big|\sum_{m\sim M} e(xm+ym^2+\lambda Nzm^3+\lambda wm^4)\Big|^p dxdydzdw
\ee
to which Theorem 1 is applicable.
In the above $\lambda$ plays the role of a parameter, nothing that $I_p(\lambda)$ decreases with $\lambda$ for $p$ an even
integer.

\subsection
{Preliminary decoupling}
\hfill\break

Denote $S=\sum_{n\sim N} e(n^2x+\lambda n^4 y)$ and $S_I=\sum_{n\in I} e(n^2x+\lambda n^4y)$ for
\hfill\break $I\subset \big[\frac N2, N\big]$ an interval. Assuming
\be\label{5.7}
\lambda N^4>\frac {N^2}{M^2}, \ \text { i.e. }  \ \lambda N^2 M^2>1
\ee
the decoupling theorem for curves gives for $p\geq 6$
\be\label{5.8}
\Vert S\Vert_p \ll N^\ve \Big(\frac NM\Big)^{\frac 12 -\frac 3p} \Big(\sum_I\Vert S_I\Vert_p^2\Big)^{\frac 12}
\ee
with $\{I\}$ a partition of $\big[\frac N2, N\big]$ in $M$-intervals.
Hence
\be\label{5.9}
I_p(\lambda)\ll N^\ve\Big(\frac NM\Big)^{p-3}\Big(\frac 1N\sum_{n\sim N}\int^1_0\int^1_0 |S_{[n, n+M]}
(x, y)|^p dxdy\Big)
\ee
where
\be\label{5.10}
|S_{[n, n+M]} (x, y)| =\Big|\sum_{m\sim M} e\big((2nx+ 4\lambda n^3y)m+(x+6\lambda n^2y)m^2+4\lambda nm^3y+\lambda m^4y\big)\Big|.
\ee

\subsection
{Distributional considerations}
\hfill\break

In view of \eqref{5.9}, \eqref{5.10} and exploiting the additional average over $n$, it is natural to analyze the distribution induced by
the map
$$
\vp: [0, 1]\times[0, 1] \times\{n\sim N\}\to \mathbb T\times\mathbb T\times [0, 4N]\times [0, 1]
$$
\be\label{5.11}
(x, y, n)\mapsto (2nx+4\lambda n^3y, x+ 6\lambda n^2y, 4ny, y)=(x', y', z', \omega').
\ee

For the time, restrict $y$ to $\big[\frac 12, 1\big]$ and denote $\mu$ the (normalized) image measure of $\vp$.
A translation $x\mapsto x-2\lambda n^2y$ $(\mod 1)$ clearly permits to replace $\vp$ by the map
$$
(x, y, n)\mapsto (2nx, x+4\lambda n^2 y, 4ny, y)
$$
and we need to analyze the distribution of $\mu$ at scale $\frac 1M\times \frac 1{M^2}\times \frac 1{\lambda M^3}
\times \frac 1{\lambda M^4}$.

Hence, let $k,\ell \in\mathbb Z, |k|\lesssim M, |\ell|\lesssim M^2$ and $\xi, \eta\in\mathbb R, |\xi|\lesssim \lambda M^3, |\eta|<\lambda M^4$.
The Fourier transform $\hat\mu$ of $\mu$ is given by
$$
\hat\mu(k, \ell, \xi, \eta) =\frac 1N \sum_{n\sim N}\int\int dxdy \, e(2nxk+ (x+4\lambda n^2y)\ell +4ny\xi+y\eta)
$$
implying
$$
|\hat\mu (k, \ell, \xi, \eta)|\ll \frac 1N\sum_{n\sim N} 1_{[2nk+\ell=0]} \, 1_{[|4\lambda n^2\ell+4n\xi+\eta|\ll N^\ve]}.
$$

It follows from the restrictions on $\xi, \eta$ that
$$
|\ell|\ll \frac {N^\ve}{\lambda N^2} +\frac {M^3}{N} =\frac {M^4}{N^2} <\frac 1{\lambda N^{2-\ve}} +\frac {M^3} N
$$
and
$$
|k|\ll\frac 1{\lambda N^{3-\ve}} +\frac {M^3}{N^2}.
$$
Assume further
\be\label{5.12}
\lambda>N^{-3+\ve} \text { and } \ M<N^{\frac 23-\ve}
\ee
as to ensure $k=\ell=0$.

Hence $\mu\ll N^\ve \pi_{z', w'}[\mu]$.
Returning to \eqref{5.9}, we may therefore bound
$$
\frac 1N\sum_{n\sim N} \int_0^1 \int_{\frac 12}^1 dxdy\Big|\sum_{m\sim M} e\big( (2nx+4\lambda n^3y)m+(x+6\lambda n^2 y)m^2+ 4\lambda nym^3+\lambda m^4y\big)|^p 
$$
by
\be\label{5.13}
\frac 1{N^{1-\ve}}\sum_{n\sim N}\int^1_0\int^1_0 \int^1_{\frac 12} dx'dy'dy\Big(\sum_{m\sim M} \, e(x'm+y'm^2+4\lambda
nym^3+\lambda m^4 y)\Big)\Big|^p
\ee

Since $[\frac 12, 1]\times\{n\sim N\}\to [0, 4N]\times [0, 1]: (y, n)\mapsto (4ny, y)$ induces a measure bounded by
the     uniform measure
at scale $1\times dw'$, it follows that at scale $\frac 1M\times \frac 1{M^2}\times \frac 1{\lambda M^3}\times\frac 1{\lambda
M^4}$ $\mu$ may be majorized by uniform measure up to a factor $N^\ve(1+\lambda M^3)$.
Hence \eqref{5.13} may be bounded by
\be\label{5.14}
N^\ve(1+\lambda M^3)\int^1_0\int^1_0\int_0^1 \int_0^1\Big|\sum_{m\sim M} e(x'm+y'm^2+\lambda Nz' m^3 +\lambda
w'm^4)\Big|^p dx'dy'dz'dw'.
\ee

One may do better.
Assume $\lambda M^3>100$ and shift in \eqref{5.13} the $y$-variable by $o\big(\frac 1{\lambda M^4}\big)$, i.e.
replace $y$ by $y+\frac z{\lambda M^4}, z= o(1)$.
One obtains
$$
\frac 1{N^{1-\ve}} \sum_{n\sim N}\int^1_0\int^1_0\int^1_{\frac 12}\int^1_0 
dx'dy'dydz \Big|\sum_{m\sim M}\, e(x'm+y'm^2+4\lambda n\Big(y+\frac z{\lambda M^4}\Big) m^3 +\lambda m^4 y)\Big|^p.
$$
Assuming
\be\label {5.15}
\lambda M^4<N
\ee
we note that for fixed $\frac 12\leq y\leq 1$, the map $(n, z)\mapsto n\big(y+\frac z{\lambda M^4}\big)$
induces a normalized measure essentially bounded by $\frac 1N 1_{[0, 2N]}$.
Consequently, under the condition \eqref{5.15}, \eqref{5.13} is bounded by
\be\label{5.16}
N^\ve\int_0^1\int^1_0\int^1_0\int^1_0 \Big|\sum_{m\sim M} e(x' m+y'm^2+\lambda Nz' m^3+\lambda w' m^4\Big|^p
dx'dy'dz'dw'.
\ee
Taking $M<N^{\frac 23}$, \eqref{5.15} will hold for $\lambda<N^{-5/3}$.

\subsection
{Application of mean value theorems}
\hfill\break

Use the 8th moment bound (\cite{B-I2}), or equivalently, Theorem 2 in the paper, we get
\be\label{5.17}
\max_{|c_m|\leq 1} \int\int\int dx'dy'dz'\Big|\sum_{m\sim M} c_m \, e\Big(x'm+y'm^2+\frac {z'}{M} m^3\Big)\Big|^8 \ll M^{4+\ve}
\ee
Application of \eqref{5.17} to \eqref{5.16} with fixed $w'$ and $p=8$ implies then
\be\label{5.18}
\int_0^1 \int^1_{\frac 12} \Big|\sum_{n\sim N} e(n^2x+\lambda n^4y)\Big|^8 dxdy \ll N^\ve \Big(\frac NM\Big)^5 M^4 
\Big( 1+\frac 1{\lambda NM^2}\Big) \ll N^{4+\frac 13+\ve}+ \frac{N^{2+\ve}} \lambda
\ee
taking $M=N^{\frac 23-\ve} $ and $ N^{-3+\ve}< \lambda < N^{-2}$.

Braking up the range $y\in [0, 1]$ in sub-intervals $\big[\frac 12\sigma, \sigma\big], \sigma =2^{-s}$ 
a change of variables and replacement of $\lambda$ by $N^{-\frac 73}\sigma$ in \eqref {5.18} gives
$$
\int^1_0\int^1_{N^{-\frac 23}} \Big| \sum_{n\sim N}  e(n^2x+N^{-\frac 73}n^4 y)\Big|^8 dxdy \ll N^{\frac {13}3+\ve}.
$$
The remaining range is simply bounded by
$$
N^{-\frac 23}I_8(N^{-3}) \leq N^{\frac 43}I_6 (N^{-3})\ll N^{\frac {13}3+\ve}.
$$
Hence we establish 0.17.

\begin{theorem}\label{Theorem12}
\be\label{5.19}
I_8 \leq I_8 (N^{-\frac 73}) \ll N^{\frac {13}3+\ve}
\ee
\end{theorem}
\medskip

Next, one may consider the 10th moment.
Setting $p=10$ in \eqref{5.9} implies with $M=N^{\frac 23-\ve}, N^{-\frac 83}<\lambda< N^{-\frac 53}$
\begin{align}\label{5.20}
&\int^1_0\int^1_{\frac 12} \Big|\sum_{n\sim N} e(n^2x+\lambda n^4y)\Big|^{10} dxdy\ll\nonumber\\
&N^\ve \Big(\frac NM\Big)^7\int^1_0\int^1_0\int^1_0\int^1_0 \Big|\sum_{m\sim M}
e(mx+m^2 y+\lambda Nm^3z+\lambda m^4w)\Big|^{10} dxdydzdw.
\end{align}

 Apply Theorem 4  with $\vp_2(t)=t^3, \vp_3(t)=t^4$ and
$\delta =\lambda^{-1}  N^{-1} M^{-3}$, $\Delta =\lambda^{-1} M^{-4}$ 

This gives the bound
$$
N^\ve \Big(\frac NM\Big)^7 \{\delta \Delta ^{3/4} M^7+(\delta+\Delta)M^6 +M^5\}.
$$
\be\label{5.21}
\ll N^\ve(N^2 \lambda^{-7/4}+ N^{\frac {11} 3} \lambda^{-1} +N^{\frac {17}3})\ll N^{4+\ve} \lambda^{-1}
\ee
for $\lambda$ as above.

Thus
\be\label{5.22}
\int^1_0\int^1_{N^{-1}} \Big|\sum_{n\sim N} e(n^2x+ N^{-\frac 53} n^4y)\Big| dxdy\ll N^{\frac {17}3+\ve}.
\ee
The remaining range may be captured using \eqref{5.19}, i.e.
$$
\int_0^1\int_0^{N^{-2/3}}\Big|\sum_{n\sim N} e(n^2x+N^{-\frac 53}n^4y)\Big|^{10} dxdy\ll
N^{-\frac 23}N^2 I_8 (N^{-\frac 73})\ll N^{\frac {17}3+\ve}.
$$
Hence we establish Theorem 11.

\end{document}